\documentclass[oneside,9pt, a5paper]{article}
\usepackage[english]{babel}
\usepackage{amsmath, amssymb, amsthm,accents}
\usepackage{graphicx}
\usepackage[top=2cm, left=1.5cm, marginparwidth=2cm]{geometry}
\usepackage{cite}
\usepackage{xcolor}
\usepackage{caption}

    \advance\textwidth  by -31truemm     
    \advance\textheight by -29truemm    

\newcommand{\bs}{\boldsymbol}
\newcommand{\bal}{\bs\alpha}
\newcommand{\bbt}{\bs\beta}
\newcommand{\bg}{\bs\gamma}
\newcommand{\be}{\bs e}
\newcommand{\ba}{\bs a}
\newcommand{\bv}{\bs v}

\newcommand{\bM}{\bs M}
\newcommand{\bF}{\bs F}

\newcommand{\x}{\times}
\newcommand{\br}{\bs r}

\newcommand{\om}{\omega}
\newcommand{\bom}{\bs\om}
\renewcommand{\d}{\partial}
\newcommand{\I}{{\bf I}}
\newcommand{\J}{{\bf J}}

\newcommand{\Q}{{\bf Q}}
\newcommand{\q}{\quad}
\newcommand{\qq}{\qquad}
\newcommand{\g}{{\rm g}}
\newcommand{\lm}{\lambda}
\newcommand{\gam}{\gamma}
\newcommand{\vfi}{\varphi}
\newcommand{\tk}{\widetilde{k}}
\renewcommand{\th}{\vartheta}
\DeclareMathOperator{\diag}{diag}
\DeclareMathOperator{\const}{const}

\DeclareMathSymbol{\widetildesym}{\mathord}{largesymbols}{"65}

\newtheorem*{pro}{Proposition}

\theoremstyle{remark}
\newtheorem*{rem}{Remark}

\begin{document}

\title{The influence of the rolling resistance model on tippe top inversion}
\date{}
\author{A.\,A.\,Kilin, E.\,N.\,Pivovarova}

\maketitle


\begin{abstract}

In this paper, we analyze the effect which the choice of a friction model has on tippe top inversion in the case
where the resulting action of all dissipative forces is described not only by the force applied at the contact point, but
also by the additional rolling resistance torque. We show that, depending on the
friction model used, the system admits different first integrals. In particular, we
give examples of friction models where the Jellett integral,
the Lagrange integral or the area integral is preserved.

We examine in detail the case where the action of all dissipative forces reduces to the horizontal rolling resistance torque.
For this case we find permanent rotations of the system and analyze their linear stability. Also, we show that for this friction
model no inversion is observed.

\end{abstract}

\section*{Introduction}

In this paper we analyze the dynamics of the tippe top on a smooth
plane in the presence of forces and torques of rolling resistance. Tippe top inversion
has attracted the attention of researchers for the last decades~\cite{tt1, tt2,tt3,tt4,karap, wojc, langerock,
langerock_prototype, brien, karap_rub1, karap_rub2,wojc2}. Reference~\cite{langerock_prototype}
discusses the development of a spherical prototype of the top which a capable of performing various modes of motion (in particular,
complete or partial inversion)
by changing the mass-geometric
characteristics. Many studies from the last century~\cite{tt1,
tt3,tt4, brien, karap_rub1, karap_rub2} and from the early part of this century~\cite{tt2,karap,wojc, langerock,wojc2}
gave mathematical explanations of tippe top inversion. They investigated the stability of steady-state
(dissipation-free) solutions, analyzed parameter values, and found cases where tippe top inversion is possible.

A detailed analysis of the dynamics of an axisymmetric top is possible
since the system has the Jellett integral. In his ``Treatise on the Theory
of Friction''~\cite{jellett} Jellett pointed out that some quantity (named
later after him) in the problem of the motion of a body of revolution on a
plane remains unchanged when adding an arbitrary friction force applied at
the point of contact. There exist various generalizations and analogs of
the Jellett integral. As an example, we mention the integral found
in~\cite{twonon} for a rigid spherical shell with internal body. 

Most studies of tippe top dynamics use the classical sliding friction model, which is proportional
to the velocity of the point of contact of the tippe top with the plane. However,
various friction models have been proposed recently
for a more accurate analysis which provides not only a qualitative, but also a quantitative description
of the system dynamics. For example, a comparative analysis of the most frequently
used friction models is made in Ref.~\cite{leine}. The authors of Ref.~\cite{rollfr}
present a model of viscous rolling friction which gives a fairly accurate description
of the rolling motion of spherical bodies on a horizontal plane. Also, this friction
model explains qualitatively some dynamical effects, in particular,
retrograde motion of a rolling disk at the final stage~\cite{ufn}.

In this paper we address the problem of the influence of the friction model
on tippe top inversion. In particular, we examine the situation where the resultant action
of all dissipative forces is described not only by the force applied
at the point of contact, but also by an additional rolling resistance torque.
It turns out that, depending on the chosen friction model,
the system admits different first integrals: the Jellett integral, the Lagrange
integral or the area integral.

In this work we also carry out a qualitative analysis of
tippe top dynamics in the case
where the action of dissipative forces reduces to the horizontal rolling resistance torque.
We show that, in this friction model, no tippe top inversion is observed.

\newpage

\section{Equations of motion and conservation laws}

\subsection{Formulation of the problem}

\begin{figure}[!ht]
\centering\includegraphics[scale=1]{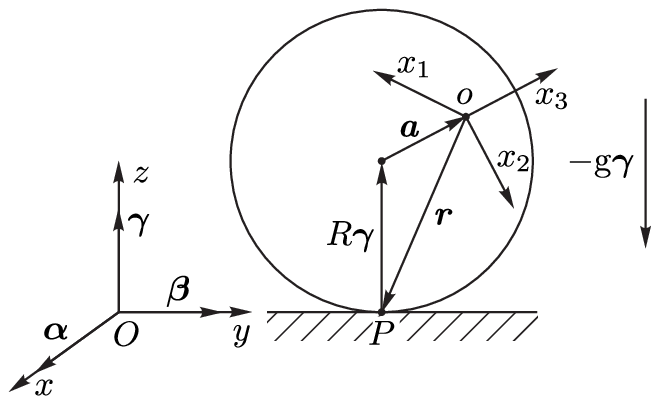}
\caption{\label{fig1}}
\end{figure}

Consider the motion of a heavy unbalanced ball of radius $R$ and mass $m$ with
axisymmetric mass distribution which rolls with slipping
on a horizontal plane under the action of gravity (Fig.~\ref{fig1}). The
system is acted upon by different resistance forces, which depend on the
type of coating of the contacting surfaces, air resistance etc.
As is well known, this system of forces generally reduces to
the resultant of resistance forces, $\bF$, and to the resistance torque $\bM_f$. We
assume that
in this case the motion of the ball is subject to the following
assumptions:\vspace{-2mm}
\begin{enumerate}\itemsep=-2pt
\item[--] the ball contacts the surface at one point $P$;
\item[--] the resultant of resistance forces, $\bF$, is applied to the point of contact;
\item[--] the ball is acted upon by the principal rolling resistance torque $\bM_f$,
which includes torque $\br\x\bF$,
    but may not be equal to it in the general case.
\end{enumerate}

To describe the motion of the ball, we introduce two coordinate systems:\vspace{-1mm}
\begin{enumerate}\itemsep=-2pt
\item[--] a fixed (inertial) coordinate system $Oxyz$ with origin on the supporting
plane and with the axis $Oz$ directed vertically upwards.
\item[--] a moving coordinate system $ox_1x_2x_3$ attached to the ball, with origin
at the center of mass of the system and with the axis $ox_3$ directed along the symmetry axis
of the ball.
\end{enumerate}
In what follows, unless otherwise specified, all vectors will be referred to the
moving axes $ox_1x_2x_3$.

We assume that the center of mass of the system is displaced relative to the
geometric center of the ball along its symmetry axis by distance $a$ and is given
by the vector $\ba=(0,0,a)$.

Let us denote the projections of the unit vectors directed along the fixed axes
$Oxyz$ onto the axes of the moving coordinate system $ox_1x_2x_3$ as
follows:
\[\bal=(\alpha_1,\alpha_2,\alpha_3), \q \bbt=(\beta_1,\beta_2,\beta_3), \q \bg=(\gam_1,\gam_2,\gam_3).\]
The orthogonal matrix $\Q\in SO(3)$ whose columns are the coordinates of the vectors
$\bal,\bbt,\bg$ specifies the orientation of the body in space.

\begin{rem}
By definition, the vectors $\bal$, $\bbt$ and $\bg$ satisfy the relations
\[\begin{gathered}(\bal,\bal)=1,\q(\bbt,\bbt)=1,\q(\bg,\bg)=1,\\ (\bal,\bbt)=0,\q(\bbt,\bg)=0,\q(\bg,\bal)=0.\end{gathered}\]
\end{rem}

Let $\bv$ be the velocity of the center of mass of the ball, and let $\bom$ be its
angular
velocity, both defined in the coordinate system $ox_1x_2x_3$. Then the evolution
of the orientation and of the position of the ball is given by the kinematic relations
\begin{equation}\begin{gathered}
\dot\bg = \bg\x\bom,\qq\dot\bal = \bal\x\bom,\qq \dot\bbt = \bbt\x\bom,\\
\dot x = (\bv,\bal),\qq \dot y = (\bv,\bbt),
\end{gathered}\label{kin_eq}
\end{equation}
where $x$ and $y$ are the coordinates of the center of mass $o$ of the ball in the
fixed coordinate system
$Oxyz$. The coordinate $z$
of the center of mass is uniquely defined from the condition that the ball
move without loss of contact with the plane
\begin{equation}z+(\br,\bg)=0,\label{constr_g}\end{equation}
where the radius vector
 of the contact point $\br$ in the axes $ox_1x_2x_3$ can be represented in
the form \begin{equation}\br=-R\bg-\ba.\label{r}\end{equation}

\newpage
\subsection{Equations of motion}

Differentiating Eq.~\eqref{constr_g} taking~\eqref{r} and the relation
$\dot z = (\bv,\bg)$ into account, we obtain a (holonomic) constraint equation in the form
\begin{equation}f=(\bv+\bom\x\br,\bg)=(\bv_p, \bg)=0,\label{constr}\end{equation}
where
 $\bv_p$ is the velocity of the ball at the point of contact with the plane.

We write the Newton\,--\,Euler equations for changing the linear and angular momenta of the ball
in the form
\begin{equation}
\begin{aligned}
& m\dot\bv + m\bom\x \bv  ={ N}\frac{\d f}{\d \bv} + \bF - m\g\bg,&\\
&\I\dot\bom + \bom\x \I\bom= {N}\frac{\d f}{\d \bom}+ \bM_f,& 
\end{aligned}\label{lagreqslm}
\end{equation}
where $\I=\diag(i_1, i_1, i_3)$ is the central tensor of inertia of the ball, $\g$
is the free-fall acceleration, $\bF$ is the friction force applied at the
point of contact, $\bM_f$ is the rolling resistance torque, $ N\frac{\d
f}{\d \bv}$ and $ N\frac{\d f}{\d \bom}$ are the force and the torque of the reaction
of the supporting plane, respectively.

Let us express $ N$ from the first
equation of~\eqref{lagreqslm} and from the derivative of the constraint
equation~\eqref{constr} with respect to time $\dfrac{d f}{d t}=0$
\begin{equation}{ N}=m(\bv,\bg)^{\bs\cdot} - (\bF,\bg) + m\g.\label{nu}\end{equation}

If the surfaces of the ball and the plane are homogeneous (but not necessarily isotropic),
then the force $\bF$ and the torque $\bM_f$ depend only on the
variables $(\bv,\,\bom,\,\bg)$. In this case, the system of equations describing the
change of
variables $(\bv,\,\bom,\,\bg)$ closes in itself and takes the form
\begin{equation}
\begin{aligned}
& m\dot\bv + \bom\x m\bv + m(\bom,\br\x\bg)^{\bs\cdot}\bg  =\bF_{\|},&\\
&\J\dot\bom + \bom\x \I\bom+ m((\bom, (\br\x\bg)^{\bs\cdot}) -\g) \br\x\bg = \bM_f - \br\x\bF_\bot,&\\
&\dot\bg + \bom\x\bg = 0,&
\end{aligned}\label{lagreqs}
\end{equation}
where $\J = \I + m(\br\x\bg)\otimes(\br\x\bg)$\footnote{Symbol $\otimes$
denotes tensor multiplication of vectors, which associates the matrix ${\bf A}$ with components
$A_{i, j}=a_ib_j$ to the pair of
vectors $\bs a,\,\bs
b$.},
$\bF_{\|}=\bF-(\bF,\bg)\bg$ is the horizontal component of force $\bF$,
and $\bF_{\bot}=(\bF,\bg)\bg$ is its vertical component.

In addition, these equations must be restricted to the submanifold given
by the constraint equation~\eqref{constr} and the geometric relation
\begin{equation}
(\bv+\bom\x\br,\bg)=0,\qq \bg^2=1.
\label{geom}
\end{equation}
Since these functions are first integrals of the system~\eqref{lagreqs}, this restriction
is satisfied uniquely.

\begin{rem}
In principle, one can consider the trajectories of the system~\eqref{lagreqs}
on other level sets of the integrals~\eqref{geom}, but they have no explicit physical interpretation.
\end{rem}

Thus, Eqs.~\eqref{kin_eq} and \eqref{lagreqs} completely describe the motion
of the tippe top on a smooth plane with friction.

\subsection{Laws of resistance and cases of existence of additional integrals}

In the case of rolling on an \textit{absolutely smooth plane} ($\bF=0$, $\bM_f=0$)
the system~\eqref{lagreqs} admits, in addition to the integrals~\eqref{geom},
three more integrals of motion:
\begin{enumerate}\itemsep=-2pt
\item[--] the energy integral
\begin{equation}
\mathcal{E}=\frac12 (\bom, \I \bom) + \frac12 m(\bv, \bv)+m\g a\gam_3,
\label{int_en}
\end{equation}
\item[--] the Lagrange integral
\begin{equation}
F=i_3\om_3,
\label{int_lagr}
\end{equation}
\item[--] the Jellett integral
\begin{equation}
G=-(\J\bom,\br).
\label{int_g}
\end{equation}
\end{enumerate}
The system also admits the area integral, which in this case is a linear combination
of the integrals~\eqref{int_lagr} and~\eqref{int_g}
\begin{equation}
C=(\J\bom,\bg) =(\I\bom,\bg) = \frac1R(G-aF).
\label{int_pl}
\end{equation}

In the presence of rolling resistance forces the energy integral~\eqref{int_en} is not preserved and its evolution is governed by
the equation
\[\dot{\mathcal{E}} = (\bF, \bv) +(\bM_f,\bom).\]

In the case of dissipative forces $\bF$ and torques $\bM_f$, the following inequality holds for all possible
values of the variables $\bv,\,\bom,\,\bg$:
\[\dot{\mathcal{E}} \leqslant 0.\]

When the force $\bF$ and the torque $\bM_f$ depend arbitrarily on phase
variables, the integrals~\eqref{int_lagr} and~\eqref{int_g} disappear as well. However, there exist several particular cases where
one of them (or a combination of both) is preserved.

\begin{pro}
The system~\eqref{lagreqs}, which describes the rolling of the ball with axisymmetric
mass distribution on the plane, admits an additional motion integral, linear in
$\bom$,
for an arbitrary rolling resistance force $\bF$
under the following restrictions on the rolling resistance torque $\bM_f$:
\begin{enumerate}\itemsep=-2pt
\item[$1^\circ.$] if $(\bM_f,\br)=0$, then the Jellett
integral is preserved
\[G=-(\J\bom,\br)=\const;\]
\item[$2^\circ.$] if $(\bM_f,\be_3)=0$, then the Lagrange integral
is preserved
\[F=i_3\om_3=\const;\]
\item[$3^\circ.$] if $(\bM_f,\bg)=0$, then the area integral is preserved
\[C=(\J\bom,\bg)=\const.\]
\end{enumerate}
\end{pro}
We recall that the rolling resistance torque $\bM_f$ includes the torque $\br\x\bF$, but in the general case
it does not coincide with it.

\begin{rem}
If any two conditions of the proposition are satisfied, all three first
integrals, $G$, $F$ and $C$, are preserved. This is due to the linear dependence of
the conditions which follows from relation~\eqref{r}.
\end{rem}

Case $1^\circ$ has been extensively studied in the literature, whereas cases $2^\circ$
and $3^\circ$ have not been mentioned previously. Below we give some examples of resistance models where the
system~\eqref{lagreqs} admits an additional integral.

$1^\circ.$ $(\bM_f,\br)=0$. The most frequently used particular case of this
situation is
\[\bM_f=\br\x\bF,\]
i.e., the resistance forces reduce to the only friction force applied at the point of contact.
In this case, one usually considers the model of viscous or dry friction, for which the friction force is given by
\begin{equation}
\bF_v=-\mu\bv_p,\qq \bF_d = -\mu N\frac{\bv_p}{|\bv_p|},
\label{frict_visc}
\end{equation}
where $\bv_p$ is the velocity of the body at the contact point, $N$ is the
reaction force of the support, and $\mu$ is the friction coefficient. The
model of viscous friction is most frequently used to explain tippe top
inversion~\cite{tt1, tt2,tt3,tt4, wojc, langerock, langerock_prototype,
brien,wojc2}, to investigate the motion of many other mechanical systems
and to explain for the dynamical effects observed in them, for example,
the rise of a spinning egg (Jellett's egg)~\cite{moffat, pwegg}.

\begin{rem}
Another example of such a friction torque which satisfies the condition that the introduced forces be dissipative
(that there be no growth of energy) is the following model:
\begin{equation}
\bF=-\mu\bom\x\br,\qq\bM_f=-\mu\br\x(\bom\x\br).
\label{mom}
\end{equation}
\end{rem}

$2^\circ.$ $(\bM_f, \be_3)=0$ --- \textit{anisotropic friction} for
which the projection of the total rolling resistance torque onto the symmetry axis
of the body, $\be_3=(0,0,1)$, is zero.

The simplest example is provided by the following model of viscous anisotropic
friction:
\begin{equation}
\bF=0,\qq\bM_f=-\mu\be_3\x(\bom\x\be_3).
\label{mom_ex}
\end{equation}
Physically this model can be achieved if the body is coated with some special
material or consists of rollers so that in one direction the body can slip freely (rotate freely about the symmetry axis), while in the other direction
resistance arises. Such friction arises in omniwheels, which is why this case
can be called the omniball.

\begin{rem}
Preliminary research has shown that there is no tippe top inversion when one uses a friction
model of the form~\eqref{mom_ex} in which the angular velocity of
rotation of the ball
relative to the symmetry axis is preserved.
\end{rem}

$3^\circ.$ $(\bM_f, \bg)=0$, i.e., the projection of the total rolling resistance torque
onto the vertical is zero.
In other words, there is no spinning resistance of the ball.

In this work we examine the simplest friction model of the following form:
\begin{equation}\bF=0,\qq\bM_f = -\mu_r\bom_\bot,\label{mod_frict}
\end{equation}
where $\bom_\bot=\bom-(\bom,\bg)\bg$ is the horizontal component of the angular
velocity and $\mu_r$ is the coefficient of rolling friction. Physically this
model describes, for example, the case of fast rolling of the ball with slipping
between two horizontal planes in the presence of dry friction.
In this case, the sum of resistance forces is zero, but the rolling resistance torque
$\bM_f$ remains nonzero.

This model is the most natural after the model~\eqref{frict_visc} and may be interesting when it
comes to analyzing the possibility of tippe top inversion and comparing the results
with those of well-known models.

In the present paper, we focus on exploring tippe top dynamics within the framework of
this resistance model and show that in this model the tippe top does not flip over under fast
rotations.

\begin{rem}
The condition $(\bM_f, \bg)=0$ is also satisfied by central resistance forces,
for example, air drag force. In the chosen notation this corresponds
to motion under the action of an arbitrary central force $\bF$ and the zero
rolling resistance torque $\bM_f=0$. Consequently, the integral $C$ in this case
is preserved as well.
However, although this case leads to preservation of the area integral,
it does not satisfy the original assumption that the force is applied
at the contact point. We have mentioned it only to illustrate the resistance model considered.
\end{rem}

In proving item~$3^\circ$ of the proposition, the condition that the tensor of inertia and the direction of displacement
of the center of mass be axisymmetric is not used.
Therefore, \textit{the area integral is also preserved
in the problem of the rolling of an unbalanced dynamically asymmetric ball
in the case where the projection
of the total rolling resistance torque onto the vertical is zero.}

In a similar manner, it is easy to show that \textit{in the case $(\bM_f, \br)=0$
the Jellett integral
is also preserved in the problem of the rolling of the Chaplygin ball}
(a balanced dynamically asymmetric ball).


\section{Reduction of the equations of motion}

We now turn to a study of the dynamics of the ball using the
friction model~\eqref{mod_frict}.

For further analysis of the equations, we introduce dimensionless variables in the
following form:
\begin{equation*}t\rightarrow\sqrt{\frac{R}{\g}}\,t,\q \br\rightarrow R\br,\q \bom\rightarrow\sqrt{\frac{\g}{R}}\,\bom,\q \I\rightarrow mR^2\I,\q \ba\rightarrow R\ba,\q \mathcal{E}\rightarrow m\g
R{\mathcal{E}}.\label{br}\vspace{-2mm}\end{equation*}
Such a change of variables is equivalent to
\[m=1,\q R=1,\q \g=1.\]

In the case under consideration the equations of motion for $\bom,\bg$ decouple
from the
complete system and have the form
\begin{equation}
\begin{aligned}
&\J\dot\bom + \bom\x \I\bom+ ((\bom, \br\x(\bg\x\bom)) -1) \br\x\bg = \mu_r\bg\x(\bom\x\bg) ,&\\
&\dot\bg + \bom\x\bg = 0.&
\end{aligned}\label{eqs_mur}
\end{equation}

Due to the existence of a pair of integrals the phase space of the
system~\eqref{eqs_mur} is foliated by four-dimensional invariant submanifolds
\[\mathcal{M}_C^4 = \{(\bom,\bg)\,\|\, (\bg,\bg)=1,\,\, (\J\bom,\bg)=C\}.\]
The analysis of the flow on the submanifolds $\mathcal{M}_C^4$ becomes simpler
since the
system~\eqref{eqs_mur} has the symmetry field
\[\hat{\bs u} = \om_1\dfrac{\d}{\d \om_2} - \om_2\dfrac{\d}{\d \om_1} + \gam_1\dfrac{\d}{\d \gam_2} - \gam_2\dfrac{\d}{\d \gam_1},\]
which defines rotations about the symmetry axis of the ball $ox_3$.
This symmetry makes it possible to perform reduction (reduce the order of
the system) on the invariant submanifolds $\mathcal{M}_C^4$.
For this we proceed as follows. Let us choose variables $\gam_3,\,K_1,\,K_2,\vfi$
which parameterize $\mathcal{M}_C^4$ so that three of them are first integrals of
the symmetry field~\cite{symm}:
\[\gam_3=\gam_3,\q K_1 = i_3\om_3,\q K_2=
\frac{1}{k}(\gam_1\om_2-\gam_2\om_1),\q \vfi = \arctan\dfrac{\gam_2}{\gam_1},\]
where $k = \sqrt{\frac{1-\gam_3^2}{i_1+a^2(1-\gam_3^2)}}$. The inverse transformation
has the form
\begin{eqnarray*}
&\gam_1=\sqrt{1-\gam_3^2}\cos\vfi,\q \gam_2=\sqrt{1-\gam_3^2}\sin\vfi,\\[1mm]
&\om_1 = \dfrac{(C-\gam_3K_1)\gam_1-i_1kK_2\gam_2}{i_1(1-\gam_3^2)},\,\,\om_2=\dfrac{(C-\gam_3K_1)\gam_2+i_1kK_2\gam_1}{i_1(1-\gam_3^2)},\,\, \om_3=\dfrac{K_1}{i_3}.
\end{eqnarray*}

The evolution of the new variables on the level set $\mathcal{M}_C^4$ is governed
by the
equations
\begin{equation}
\begin{gathered}
\q \dot\gam_3=kK_2,\q \dot K_1 = -\frac{\mu_r}{i_1}\left(K_1\tk - \gamma_3C\right),\\
\dot K_2 = -\frac{k(C-\gam_3K_1)(C\gam_3 - K_1)}{i_1(1-\gam_3^2)^2}-ka - \dfrac{\mu_r K_2 k^2}{1-\gam_3^2},
\end{gathered}\label{syst_kgamma}
\end{equation}\vspace{-3mm}
\begin{equation}
\dot\vfi = \dfrac{\dot\gam_1\gam_2 - \dot\gam_2\gam_1}{\gam_1^2+\gam_2^2}=\dfrac{K_1}{i_3}-\dfrac{\gam_3(C-K_1\gam_3)}{i_1(1-\gam_3^2)},
\end{equation}
where we have introduced the notation $\tk = \left(i_1-(i_1-i_3)\bg_3^2\right) / i_3$.
Due to the special choice of variables the first three equations decouple and
form a closed reduced system. It is this system that we will investigate in what follows.

\section{Permanent rotations and their stability}

Permanent rotations of the system under consideration correspond to motion of the
ball
with a constant inclination angle of the symmetry axis relative to the vertical, i.e.,
when $\gam_3=\const$. It is obvious that
during such motions there should be no dissipation: $\dot{\mathcal{E}}=0$.
Permanent rotations are given by the equations
\begin{equation}
\dot\gamma_3=0,\q \ddot\gam_3 = 0. 
\label{perm}
\end{equation}

Solving Eqs.~\eqref{perm}, we find the following partial solutions to the
system:
\begin{enumerate}
\item Two one-parameter families of fixed points
\begin{equation}
\begin{array}{ll}
\sigma_{u\phantom{l}}:& \gam_3=1,\q K_1 = C,\q K_2 =0,\\[2mm]
\sigma_{l\phantom{u}} :& \gam_3=-1,\q K_1 = -C,\q K_2 =0.
\end{array}\label{partsol_vert}
\end{equation}
These families correspond to vertical rotations of the ball where the
center of mass is above the geometric center ($\sigma_u$) or below it
($\sigma_l$).

\item A one-parameter family of periodic solutions
\begin{equation}
\sigma_0:\q \gam_3 = -\frac{a}{c_1^2(i_1-i_3)},\q K_1 =  \frac{ai_3}{c_1(i_1-i_3)},\q K_2=0,
\label{partsol_mid}
\end{equation}
where $c_1\in(-\infty, -c_0)\cup(c_0,\infty)$, $c_0=\sqrt{a/|i_1-i_3|}$, is the
parameter of the family.
This family of periodic solutions corresponds to rotations of the ball where its symmetry axis
deviates from the vertical by the angle $\th = \arccos \gam_3$.
\end{enumerate}

The initial variables $\bom,\bg$ and the value of the integral $C$
are parameterized through $c_1$ as follows:
\[\om_1=c_1p\cos\vfi,\q\om_2=c_1p\sin\vfi,\q\om_3 = \frac{a}{c_1(i_1-i_3)}\]
\[\gam_1 = -p\cos\vfi,\q\gam_2=-p\sin\vfi,\q\gam_3=-\frac{a}{c_1^2(i_1-i_3)},\]
\[C = -c_1i_1+\frac{a^2}{c_1^3(i_1-i_3)},\]
where $p=\sqrt{1-\gam_3^2}$.

Using this parameterization and the obvious inequality
$|\gam_3|\leqslant1$, we can define the critical value of the integral $C$
at which permanent rotations arise or disappear:
\[C^* = C|_{c_1=c_0}=\dfrac{i_3\sqrt{a}}{\sqrt{|i_1-i_3|}}.\]

\subsection{Stability analysis}

Let us analyze the linear stability of partial solutions $\sigma_u$ and $\,\sigma_l$.
These solutions are families of fixed points of the complete system~\eqref{eqs_mur}.
To investigate their stability, we represent the system of differential
equations~\eqref{eqs_mur} as
\[\dot{\bs q}={\bf f}_q(\bs q),\]
where $\bs q = (\om_1,\om_2,\om_3,\gam_1,\gam_2,\gam_3)$, ${\bf f}_q(\bs q)$
is the vector whose components are functions of $\bs q$. Let us linearize
the system~\eqref{eqs_mur} near the solutions $\sigma_u,\,\sigma_l$
\[\dot{\widetilde{\bs q}}= {\bf L}_q\widetilde{\bs q},\qq {\bf L}_q=\left.\frac{\d {\bf f}_q(\bs q)}{\d \bs q}\right|_{\bs q=\bs q_{u, l}},\]
where $\widetilde{\bs q} = \bs q - \bs q_{u, l}$, and $\bs q_u$ and $\bs q_l$
are the partial solutions $\sigma_u$ and $\sigma_l$ of the
system~\eqref{eqs_mur}, respectively.

The characteristic equation of the linearized system
\[\det({\bf L}_q-\lm {\bf E})=0\]
with eigenvalues $\lm$ ($\bf E$ being an $6\x6$ identity matrix) is an equation of
degree $6$ in $\lm$ of the form
\[P_6(\lm) =\lm^2P_4(\lm)= \lm^2(a_0\lm^4+a_1\lm^3+a_2\lm^2 +a_3\lm +a_4) = 0.\]
Two zero eigenvalues of the linearized system correspond
to the geometric integral $\bg^2=1$ and to the area integral $C$ (which is
the parameter of the family).

To investigate the stability problem, we use the Routh\,--\,Hurwitz criterion for
definition of the sign of the real part of the roots of algebraic equations.

As is well known~\cite{corn}, the real parts of all roots of the equation are
negative in the case where all diagonal minors of the Hurwitz matrix are
positive
\[
\Delta_1 = a_1,\q \Delta_2=a_1a_2-a_0a_3,\q \Delta_3 = a_3\Delta_2 - a_1^2a_4,\q \Delta_4 = a_4\Delta_3
\]
under the condition that the coefficient with the highest degree, $a_0>0$, is positive.

1. \textit{The partial solution $\sigma_l$ corresponding to the lower vertical rotation.}
In the variables $\bs q$ this solution is parameterized as follows:
\[\bs q_l\,:\q\gam_1=0,\q\gam_2=0,\q\gam_3=-1,\q\om_1=0,\q\om_2=0,\q\om_3=\dfrac{C}{i_3}.\]
The coefficient with the highest degree and the diagonal minors of the Hurwitz matrix which
corresponds
to the polynomial $P_4(\lm)$ have the form
\begin{eqnarray*}
a_0&=& i_1^2,\\
\Delta_1&=& 2i_1\mu_r,\\
\Delta_2&=& 2i_1\mu_r\left(\dfrac{C^2}{i_3^2}(3i_1^2-3i_1i_3+i_3^2) +ai_1 +\mu_r^2\right),\\
\Delta_3&=& 4i_1\mu_r^2\left(\dfrac{C^2}{i_3^2}(2i_1-i_3)^2 +\mu_r^2\right)\left(a - \dfrac{C^2}{i_3^2}(i_1-i_3)\right),\\
\Delta_4&=& 4i_1\mu_r^2\left(\dfrac{C^2}{i_3^2}(2i_1-i_3)^2 +\mu_r^2\right)\left(a - \dfrac{C^2}{i_3^2}(i_1-i_3)\right)^{\!\!3}.
\end{eqnarray*}
It is easy to see that the values of $a_0,\,\Delta_1,\,\Delta_2$
are positive when $\mu_r>0$, and the values of $\Delta_3$ and $\Delta_4$ are positive
under the condition
\[
\begin{array}{l}
|C|<C^*\q \text{for}\q i_1-i_3>0,\\
|C|>-C^*\q \text{for}\q i_1-i_3<0.
\end{array}
\]

Thus, when $i_1<i_3$, the lower rotation $\sigma_l$ is always stable,
and when $i_1>i_3$, it is stable only if the area integral has values $|C|<C^*$.

2. \textit{The partial solution $\sigma_u$ corresponding to the upper vertical rotation.}
In the variables $\bs q$ this solution is parameterized as follows:
\[\bs q_u\,:\q\gam_1=0,\q\gam_2=0,\q\gam_3=1,\q\om_1=0,\q\om_2=0,\q\om_3=\dfrac{C}{i_3}.\]
The coefficient with the highest degree and the diagonal minors of the Hurwitz matrix which
corresponds to the
polynomial $P_4(\lm)$ have the form
\begin{eqnarray*}
a_0&=& i_1^2,\\
\Delta_1&=& 2i_1\mu_r,\\
\Delta_2&=& 2i_1\mu_r\left(\dfrac{C^2}{i_3^2}(3i_1^2-3i_1i_3+i_3^2) -ai_1 +\mu_r^2\right),\\
\Delta_3&=& -4i_1\mu_r^2\left(\dfrac{C^2}{i_3^2}(2i_1-i_3)^2 +\mu_r^2\right)\left(a - \dfrac{C^2}{i_3^2}(i_3-i_1)\right),\\
\Delta_4&=& -4i_1\mu_r^2\left(\dfrac{C^2}{i_3^2}(2i_1-i_3)^2 +\mu_r^2\right)\left(a - \dfrac{C^2}{i_3^2}(i_3-i_1)\right)^{\!\!3}.
\end{eqnarray*}
The values of $a_0,\,\Delta_1,\,\Delta_2$
are always positive when $\mu_r>0$, and the values of $\Delta_3$ and $\Delta_4$
are positive under the condition
\[
\begin{array}{l}
|C|>C^*\q \text{for}\q i_3-i_1>0,\\
|C|<-C^*\q \text{for}\q i_3-i_1<0.
\end{array}
\]
Consequently, when $i_3<i_1$, the upper rotation $\sigma_u$ is always unstable,
and when $i_3>i_1$, it is stable only for sufficiently large values of the
area integral, $|C|>C^*$.

3. \textit{The partial solution $\sigma_0$ corresponding to permanent rotation.}
We now investigate the stability of permanent rotations $\sigma_0$. These
rotations are periodic solutions of the complete system~\eqref{eqs_mur} and
correspond to fixed points of the reduced system~\eqref{syst_kgamma}.

To analyze the orbital stability of rotations $\sigma_0$, we examine the reduced
system~\eqref{syst_kgamma}.
To do so, we represent the system of differential
equations~\eqref{syst_kgamma} as
\begin{equation}\dot{\bs\xi} = {\bf f}_\xi(\bs\xi),\label{deq}\end{equation}
where $\bs\xi=(\gam_3,K_1,K_2)$, ${\bf f}_\xi(\bs\xi)$ is the vector whose components
are functions of $\bs\xi$. Let $\bs\xi_0$ be a partial solution
to~\eqref{deq} that corresponds to permanent rotations $\sigma_0$.

We linearize the system~\eqref{deq} near the partial solution $\bs\xi_0$ and obtain
a system of the form
\[\dot{\widetilde{\bs\xi}}= {\bf L}_\xi\widetilde{\bs\xi},\qq {\bf L}_\xi=\left.\frac{\d {\bf f}_\xi(\bs\xi)}{\d \bs\xi}\right|_{\bs\xi=\bs\xi_0},\]
where $\widetilde{\bs\xi} = \bs\xi - \bs\xi_0$, ${\bf L}_\xi$ is a linearization matrix.

The characteristic equation of the linearized system
\[\det({\bf L}_\xi-\lm {\bf E})=0\]
with eigenvalues $\lm$ is an equation of degree $3$ in
$\lm$
\[P_3(\lm) = a_0\lm^3+a_1\lm^2+a_2\lm +a_3 = 0.\]

The coefficients of the characteristic equation have the form
\begin{eqnarray*}
a_0&=& - i_1i_3 c_2(i_1-i_3)(i_1+a^2(1-c_2^2)),\\
a_1&=& -\mu_r c_2(i_1-i_3)\left(i_1i_3(1+c_2^2) +(i_1^2+i_1a^2-a^2c_2^2(i_1-i_3))(1-c_2^2)\right),\\
a_2&=& -\mu_r^2c_2(i_1-i_3)(i_1-c_2^2(i_1-i_3)) + ai_3(i_1^2 + c_2^2(i_1-i_3)(3i_1+i_3)),\\
a_3&=& a\mu_r(i_1-i_3)(1-c_2^2)(i_1+3c_2^2(i_1-i_3)),
\end{eqnarray*}
where $c_2=-\dfrac{a}{c_1^2(i_1-i_3)}$ is the value of $\gam_3$ for the partial
solution
$\sigma_0$~\eqref{partsol_mid}.

In accordance with the Routh\,--\,Hurwitz criterion, the conditions
for negativeness of the real parts of the roots of the characteristic equation
have the form
\[a_0>0,\q a_1>0,\q a_1a_2-a_0a_3>0,\q a_3>0.\]
It is easy to show that all inequalities are satisfied under the condition
$i_1>i_3$.

Thus, permanent rotations $\sigma_0$ are stable only for $i_1>i_3$ over the entire
interval of their existence.

\begin{figure}[!ht]
\centering\includegraphics[scale=0.85]{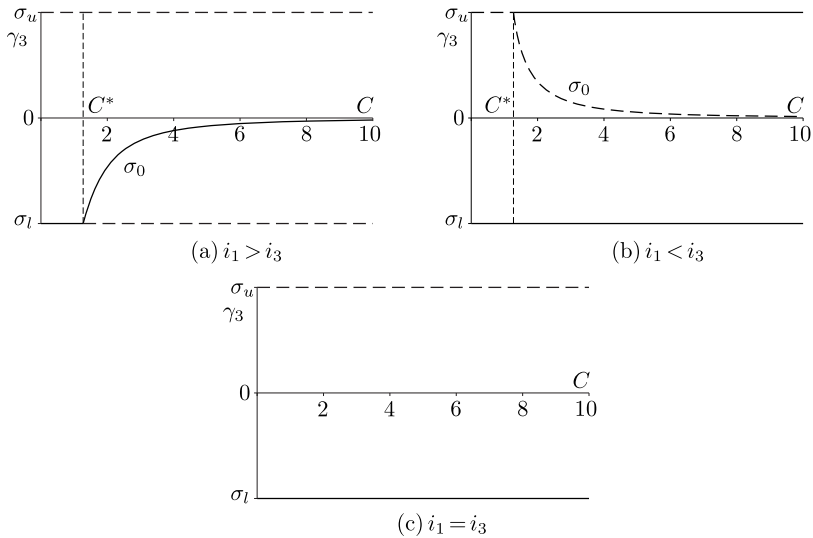}
\caption{\label{g3(k)all}Stable and unstable rotations of the tippe top versus
$C$  for different ratios of $i_1$ and $i_3$ for the parameters
$a=0.29$, $i_3=0.51$, $\mu_r=1$ and (a)~$i_1=0.55$, (b)~$i_1=0.46$,
(c)~$i_1=0.51$.}
\end{figure}

The result of linear stability analysis of the solution~\eqref{partsol_vert},
\eqref{partsol_mid} for different values of $i_1$ and $i_3$ is presented in Fig.~\ref{g3(k)all}. This figure shows
families of periodic
solutions $\sigma_u,\,\sigma_l$ and $\sigma_0$ on the plane $(C,\gam_3)$.
The solid lines represent stable permanent solutions, and the dashed lines are unstable ones. As is evident from the analysis
of eigenvalues and from the figure, there are three different cases:
\begin{enumerate}
\item $i_1>i_3$ (Fig.~\ref{g3(k)all}a). In this case the upper vertical rotation $\sigma_u$
is always unstable. The lower vertical rotation $\sigma_l$ is stable for small absolute values of the integral, $|C|<C^*$.
Permanent rotations $\sigma_0$ exist for $|C|>C^*$ and are always stable.
\item $i_1<i_3$ (Fig.~\ref{g3(k)all}b). In this case the lower vertical rotation $\sigma_l$
is always stable. The upper vertical rotation $\sigma_u$ is stable for sufficiently large absolute values of the integral,
$|C|>C^*$.
Permanent rotations $\sigma_0$ exist for $|C|>C^*$ and are always unstable.
\item $i_3=i_1$ (Fig.~\ref{g3(k)all}c). In this case the upper vertical rotation $\sigma_u$
is always unstable, and the lower vertical rotation $\sigma_l$ is always stable.
In this case there exist no permanent rotations.
\end{enumerate}

As can be seen from Fig.~\ref{g3(k)all}, the situation where on the fixed level set of the integral $C$ the lower vertical
rotation is unstable and the upper one is stable is not observed in this model.
Thus, in the friction model~\eqref{mod_frict}, a complete tippe top inversion is impossible under any initial
conditions. As we will see below, this is also confirmed by the analysis of the dependence of the energy of
permanent rotations on the area integral.

Nevertheless, the case $i_1<i_3$ admits a partial tippe top inversion when the tippe top tends to permanent rotation at a
constant inclination angle of the axis, with a small deviation from the lower vertical rotation.
However, in this case the center of mass of the tippe top always lies below the center of the ball. As
the initial energy of the tippe top increases, the critical inclination angle of the tippe top tends to $\dfrac{\pi}{2}$ $(\gam_3\rightarrow0)$.

\newpage
\section{Qualitative dynamics analysis}

A global qualitative analysis of the dynamics of the system can be carried out, for example, using
a modified Routh theory by analogy with the study of the classical model of the tippe top~\cite{karap}.
The above partial solutions which correspond to permanent and
vertical rotations $\sigma_0$, $\sigma_u$ and $\sigma_l$ can be represented in the generalized Smale diagram on
the plane $(C^2, \mathcal{E})$, where
the values of $\mathcal{E}$ correspond to the magnitude of the initial energy of the system
for a given value of the integral $C$. Fixing the level set of the integral $C$ in the diagram and
defining the initial level of energy, we can keep track of the dynamics of the system under energy dissipation.

\begin{figure}[!ht]
\centering\includegraphics[scale=0.8]{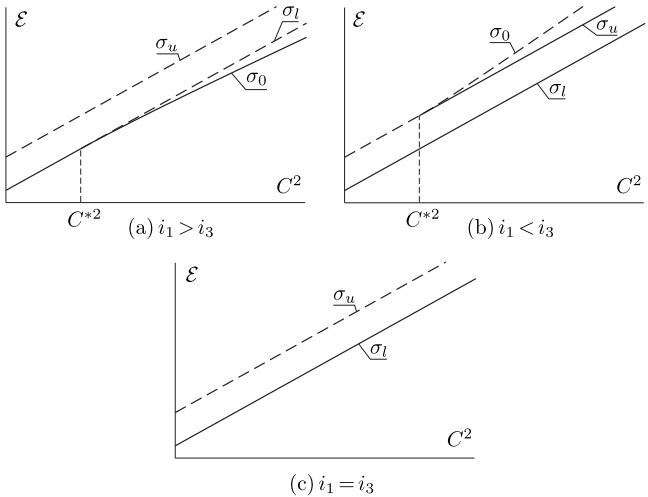}
\caption{\label{smeil}Generalized Smale diagrams
of the system~\eqref{syst_kgamma} for $a=1$, $i_3=0.5$ and (a)~$i_1=0.6$,
(b)~$i_1=0.4$, (c)~$i_1=0.5$.}
\end{figure}

Figure~\ref{smeil} shows generalized Smale diagrams corresponding to different ratios of the moments of inertia of the tippe top.
The solid lines correspond to stable steady motions, and the dashed lines to unstable ones.
On the fixed level set of the integral $C$ all
trajectories
of the system tend to stable solutions due to dissipation.

It follows from the analysis of the diagram that, depending on the system parameters and initial
conditions, the following behavior of the tippe top can be observed:
\begin{enumerate}
\item $i_1>i_3$ (Fig.~\ref{smeil}a). When the value of the integral is $|C|<C^*$,
almost all
trajectories tend to lower vertical rotations $\sigma_l$,
and when $|C|>C^*$, they tend to permanent rotations $\sigma_0$.
\item $i_1<i_3$ (Fig.~\ref{smeil}b). When the value of the integral is
$|C|<C^*$,
    almost all trajectories tend to lower vertical rotations
    $\sigma_l$. When $|C|>C^*$, the trajectories tend either to
    lower ($\sigma_l$) or to upper ($\sigma_u$) vertical
    rotations. We note that the investigation of the domains of attraction of either of the stable solutions is
a topic in its own right and goes beyond the scope of this paper. This problem can be solved, for example, by
constructing charts of dynamical regimes on the plane of initial conditions.
\item $i_3=i_1$ (Fig.~\ref{smeil}c). In this case, almost all
    trajectories tend to lower vertical rotations
    $\sigma_l$.
\end{enumerate}

Thus, no tippe top inversion is possible
for the resistance model considered.

\begin{figure}[!ht]
\centering
\begin{minipage}{.48\textwidth}
  \centering
  \includegraphics[scale=0.75]{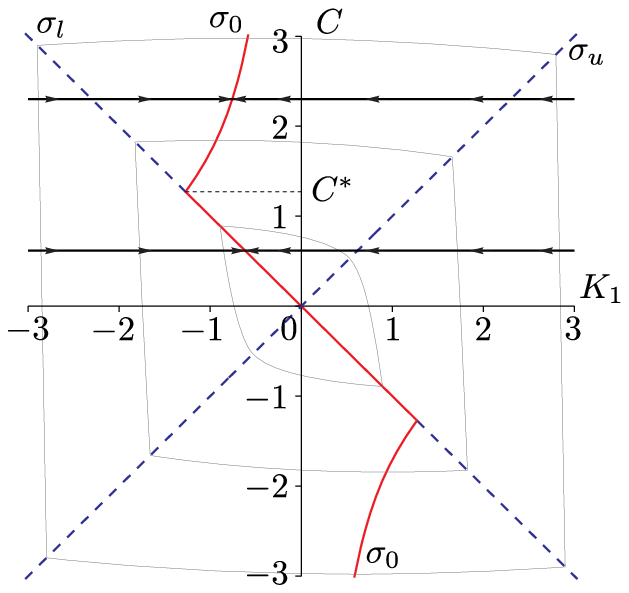}
  \captionof{figure}{Projections of the phase trajectories of the system~\eqref{syst_kgamma} onto the plane $(K_1,
C)$ for the case $i_1>i_3$ with the parameters $a=0.29$, $i_1=0.55$, $i_3=0.51$, $\mu_r=1$.}
  \label{k1k2all}
\end{minipage}\hfill
\begin{minipage}{.48\textwidth}
  \centering
  \includegraphics[scale=0.75]{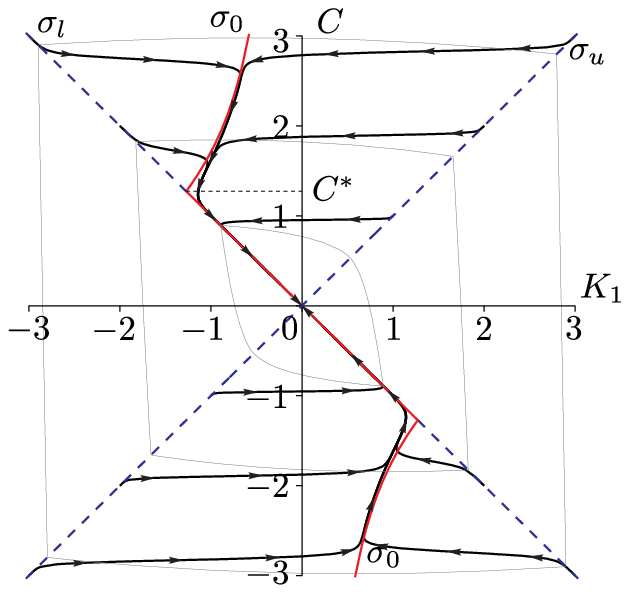}
  \captionof{figure}{Projections of the phase trajectories of the system~\eqref{eqs_mus} onto the plane $(K_1,
C)$ for the case $i_1>i_3$ with the parameters $a=0.29$, $i_1=0.55$, $i_3=0.51$, $\mu_r=1$, $\mu_s=0.001$.}
  \label{mus}
\end{minipage}
\end{figure}

We illustrate the dynamics of the system by constructing a projection
of the phase flow onto the plane $(K_1, C)$ for the case $i_1>i_3$, when a partial
tippe top inversion is observed (Fig.~\ref{k1k2all}).

As above, the solid lines indicate stable solutions and the dashed lines
represent unstable ones. The thin solid lines correspond to levels of constant energy $\mathcal{E}=\const$. The arrows denote the directions
of the system trajectories.

As can be seen from the figure, on the fixed level set of the integral $C$ all trajectories
tend to stable solutions (permanent $\sigma_0$ or lower
vertical $\sigma_l$ rotations).

\subsection{Motion with spinning friction torque}

 To conclude, we briefly examine the case where the additional integral $C$
 ceases to exist by incorporating into the system spinning friction torque in the form
\[\bM_s = -\mu_s(\bom,\bg)\bg,\]
where $\mu_s$ is the coefficient of spinning friction. We assume that the spinning
friction is much smaller than the rolling friction, i.e., $\mu_s\ll\mu_r$.

The equations of motion of~\eqref{syst_kgamma} involving rolling and
spinning friction torques take the form
\begin{equation}
\begin{aligned}
& \dot\gam_3=kK_3,&\\
& \dot K_1 = -\frac{\mu_r}{i_1}\left(K_1\tk - \gamma_3C\right) - \dfrac{\mu_s\gam_3}{i_1i_3}\left(K_1\gam_3(i_1-i_3)+Ci_3\right),&\\
&\dot K_2 = -\frac{k(C-\gam_3K_1)(C\gam_3 - K_1)}{i_1(1-\gam_3^2)^2}-ka - \dfrac{\mu_r K_2 k^2}{1-\gam_3^2},&\\
& \dot C = - \dfrac{\mu_s}{i_1i_3}\left(K_1\gam_3(i_1-i_3)+Ci_3\right).&
\end{aligned}\label{eqs_mus}
\end{equation}

As an example, we consider the effect of adding small spinning friction on the global
dynamics of the system in the case of a partial tippe top inversion.

In this case (under the condition $\mu_s\ll\mu_r$), on times $t<T_0$, the area
integral $C$ can be taken to be approximately constant. The characteristic time
$T_0$ is determined by the friction coefficient $\mu_r$. The dynamics on times
$t<T_0$ is the tendency of the tippe top to permanent
rotations $\sigma_0$. In this case, a partial inversion (elevation of the symmetry
axis) of the tippe top is observed.

On times $t\gg T_0$, an adiabatic change occurs in the value of the
integral $C$. The dynamics of the tippe top is characterized by slow motions of
the system
first along the families of permanent
rotations $\sigma_0$ (when
$|C|>C^*$) and then along the lower vertical rotations (when $|C|<C^*$). The
tippe top comes to a stop at the lower equilibrium point.

As an illustration, we present projections of the phase trajectories
of the system~\eqref{eqs_mus} onto the plane $(K_1, C)$ (Fig.~\ref{mus}). As is seen from the
figure, all trajectories tend to
permanent rotations $\sigma_0$ at initial $|C|>C^*$, then run along this solution and,
when $|C|<C^*$, evolve onto the solution corresponding to the lower vertical
rotation $\sigma_l$. At initial $|C|<C^*$ the trajectories immediately tend to
the solution $\sigma_l$.

\section*{Conclusion}
In this paper we have shown that, depending on the model
of rolling resistance, the tippe top can exhibit different behavior when it rolls
on a plane.

Within the framework of the friction
model~\eqref{mod_frict},
the area integral is preserved and no tippe top inversion is possible.
Preliminary research has shown that in the friction model~\eqref{mom_ex},
in which the Lagrange integral is preserved, no
inversion is observed either.
Interestingly, when the Jellett integral (which is a linear combination
of the Lagrange integral and the area integral) is preserved, inversion does occur.

Further research can be aimed at a more mathematically
rigorous formulation and proof of the properties of the global dynamics of the system.


\end{document}